\documentclass[a4paper,12pt,leqno]{article}
\usepackage{amsmath,amsfonts,amsthm,amssymb,dsfont}
\usepackage[alphabetic]{amsrefs}
\usepackage[OT4]{fontenc}

\long\def\symbolfootnote[#1]#2{\begingroup%
\def\thefootnote{\fnsymbol{footnote}}\footnote[#1]{#2}\endgroup}

\newtheorem{theorem}{Theorem}[section]
\newtheorem{lemma}[theorem]{Lemma}
\newtheorem{sublemma}[theorem]{Sublemma}

\newtheorem{prop}[theorem]{Proposition}
\newtheorem{cor}[theorem]{Corollary}
\newtheorem{quest}[theorem]{Question}
\theoremstyle{definition}
\newtheorem{rem}[theorem]{Remark}
\newtheorem{defin}[theorem]{Definition}
\renewcommand{\proof}{\medskip\par\noindent\textbf{Proof.} \ignorespaces}

\renewcommand{\qed}{\quad\hskip0pt\null\hfill$\square$\par}
\newcommand{\e}{\varepsilon}
\newcommand {\bbR} {{\mathbb{R}}}
\newcommand{\red}{\mathrm{red}}
\newcommand{\cyc}{\mathrm{cyc}}
\newcommand{\Out} {{\mathrm{Out}}}
\newcommand{\Hom} {{\mathrm{Hom}}}
\newcommand{\actson}{\curvearrowright}

\begin{document}

\begin{center}
\large\bfseries No-splitting property and boundaries of random groups
\end{center}

\begin{center}\bf
Fran\c cois Dahmani$^a$\footnote[1]{Partially supported by ANR grant ANR--06-JCJC-0099-01.},
 Vincent Guirardel$^a$\footnotemark[1] \&
Piotr Przytycki$^b$\footnote[2]{Partially supported by MNiSW
grant N201 012 32/0718, the Foundation for Polish Science, and ANR grant ZR58.
}
\end{center}

\begin{center}\it
$^a$
Institut de Math\'ematiques de Toulouse,
Universit\'e de Toulouse et CNRS (UMR 5219),
118 route de Narbonne,
31062 Toulouse cedex 9,
France,
\emph{e-mail:}\texttt{francois.dahmani@math.univ-toulouse.fr}
\emph{e-mail:}\texttt{vincent.guirardel@math.univ-toulouse.fr}
\end{center}

\begin{center}\it
$^b$ Institute of Mathematics, Polish Academy of Sciences,
 \'Sniadeckich 8, 00-956 Warsaw, Poland
 \\
 IRMA, CNRS(UMR 7501)

\emph{e-mail:}\texttt{pprzytyc@mimuw.edu.pl}
\end{center}

\begin{abstract}
\noindent
We prove that random groups in the Gromov density model, at any density, satisfy property $\mathrm{(FA)}$, i.e. they do not act non-trivially on trees. This implies that their Gromov boundaries, defined at density less than $\frac{1}{2}$, are Menger curves.
\end{abstract}



\section{Introduction}

The density model for random groups was introduced by Gromov. We adopt the following language from a survey by Ollivier.

\begin{defin}[{\cite[Section 9.B]{Gro93}, \cite[Definition 7]{O}}]
\label{model}
Let $F_n$ be the free group on $n\geq 2$ generators $s_1, \ldots, s_n$. For any integer $L$ let $R_L\subset F_n$ be the set of reduced words of length $L$ in these generators.

Let $0<d<1$. A \emph{random set of relators at density $d$, at length $L$} is a $\lfloor(2n-1)^{dL}\rfloor$-- tuple of elements of $R_L$, randomly picked among all elements of $R_L$.

A \emph{random group at density $d$, at length $L$} is the group $G$ presented by $\langle S |R\rangle$, where $S=\{s_1,\ldots,s_n\}$ and $R$ is a random set of relators at density $d$, at length $L$.

Let  $I\subset \mathbb{N}_+$. We say that a property of $R$, or of $G$, occurs \emph{with $I$--overwhelming probability} (shortly, \emph{w.$I$--o.p.}) \emph{at density $d$} if its probability of occurrence tends to $1$ as $L\rightarrow \infty$, for $L\in I$ and fixed $d$. We omit writing "$I$--" if $I=\mathbb{N}_+$.
\end{defin}

Note that the relators in $R_L$ need not be cyclically reduced. The case of another model is discussed in Section \ref{cyclically reduced}.

\medskip
Gromov proved the following.

\begin{theorem}[{\cite[Section 9.B]{Gro93}, \cite[Theorem 1]{Oll04}}]
\label{Gromov}
A random group is with overwhelming probability
\begin{enumerate}
\item[(i)] trivial or $\mathbb{Z}/2\mathbb{Z}$ at density greater than $\frac{1}{2}$,
\item[(ii)] word-hyperbolic, with aspherical presentation complex, at density less than $\frac{1}{2}$.
\end{enumerate}
\end{theorem}

Consequently (see e.g. \cite[Section I.3.b]{O}) with overwhelming probability at density less than $\frac{1}{2}$
a random group is torsion free, of cohomological dimension $2$, and its Euler characteristic is positive.

\medskip
We address the following question. At density less than $\frac{1}{2}$, what is the boundary at infinity of a random group $G$?

\medskip
Since $G$ is $2$--dimensional, its boundary has topological dimension $1$ (by \cite[Corollary 1.4]{BM}). The list of possibilities for the boundary is therefore limited in view of the following.

\begin{theorem}[{\cite[Theorem 4]{KK}}]
\label{KAKA}
Let $G$ be a hyperbolic group which does not split over a finite or virtually cyclic subgroup, and suppose $\partial_{\infty}G$ is $1$--dimensional. Then one of the following holds:
\begin{enumerate}
\item[(1)] $\partial_{\infty}G$ is the Menger curve;
\item[(2)] $\partial_{\infty}G$ is the Sierpi\'nski carpet;
\item[(3)] $\partial_{\infty}G$ is $\mathbb{S}^1$ and $G$ maps onto a Schwartz triangle group with finite kernel.
\end{enumerate}
\end{theorem}

Moreover, Kapovich and Kleiner prove (see \cite[Theorem 5(5)]{KK}) that $\partial_{\infty}G$ is a Sierpi\'nski carpet only if $(G; H_1, \ldots, H_k)$ is a 3--dimensional Poincar\'e duality pair, where $H_i$ are stabilizers of peripheral circles of the Sierpi\'nski carpet. But since the groups $H_i$ are virtually Fuchsian \cite[Theorem 5(2)]{KK}, this implies that the Euler characteristic of $G$ is negative. Hence w.o.p. this is not the case for a random group $G$ at density less than $\frac{1}{2}$. Case (3) is also excluded, since $G$ is w.o.p. torsion free.

\medskip
In fact, at density $d<\frac{1}{24}$, it is known that w.o.p. the boundary of a random group is the Menger curve. Namely, w.o.p. at density $d<\frac{1}{24}$, a random group satisfies  $C'(\frac{1}{12})$ small cancellation condition (see \cite[Section 9.B]{Gro93}). Champetier's theorem \cite[Theorem 4.18]{Ch} states that this condition, together with the property that each word of length $12$ is contained as a subword in one of the relators (this holds w.o.p. for random groups at any density), implies that the boundary is the Menger curve. But $C'(\frac{1}{12})$ small cancellation condition fails w.o.p. for a random group at density $d>\frac{1}{24}$ (see \cite[Section 9.B]{Gro93}).

\medskip
Nevertheless, from \.Zuk's theorem \cite[Theorem 4]{Z} it follows (see \cite[Section I.3.g]{O}) that a random group $G$ at density greater than $\frac{1}{3}$ satisfies w.o.p. Kazhdan's property (T). In particular w.o.p. $G$ does not split and by Theorem \ref{KAKA} its boundary is the Menger curve.

\medskip
We prove that this is the case at any density.

\begin{theorem}
\label{Menger}
Let $0<d<\frac{1}{2}$. Then with overwhelming probability, the boundary of a random group at density $d$ is the Menger curve.
\end{theorem}

Theorem \ref{Menger} is a consequence of Theorem \ref{KAKA} together with the discussion after it, and the following, which is the main theorem of the article.

\begin{theorem}
\label{heart}
Let $0<d<1$. Then with overwhelming probability, a random group at density $d$ satisfies property $\mathrm{(FA)}$.
\end{theorem}

Recall \cite[Section I.6.1]{S} that a group $G$ satisfies \emph{property $\mathrm{(FA)}$} if each action of $G$ on a simplicial tree has a global fixed point. When $G$ is finitely generated, it satisfies property $\mathrm{(FA)}$ if and only if it does not admit an epimorphism onto $\mathbf{Z}$ and does not split as a free product with amalgamation (see \cite[Chapter I, Theorem $15$]{S}).

\medskip
Here are additional corollaries to Theorem \ref{heart}.

\begin{cor}
Let $G$ be a random group at any density $0<d<1$. Then with overwhelming probability,
we have the following.
\begin{enumerate}
\item[(1)] $\Out(G)$ is finite.
\item[(2)] For any torsion free hyperbolic group $\Gamma$, $\Hom(G;\Gamma)$ is finite up to conjugacy.
\end{enumerate}
\end{cor}

Assertion (2) is equivalent to the fact that a random system of
equations at density $d$ has w.o.p. only
finitely many conjugacy classes of solutions in any torsion-free
hyperbolic group.

Assertion (1) is a well known (by experts) application of
Bestvina-Paulin argument
\cite{Be_degenerations,Pau_arboreal}
and Rips theory \cite{BF_stable,GLP1}.
Assertion (2) is stronger.
It follows from Sela's theory
\cite{Sela_diophantine7} and the fact that property (FA) is inherited by
quotients.
More precisely, if $\Gamma$ is hyperbolic, and $\Hom(G;\Gamma)$ is
infinite modulo conjugacy, then Bestvina-Paulin argument
provides an action of
$G$ on an $\bbR$-tree $T$.
This action factors through a  group $L$ (so called $\Gamma$-limit,
possibly not
finitely presented), such that $L\actson T$ is so-called superstable
(see \cite[Lemma 1.3]{Sela_diophantine7}).
By \cite{Sela_acylindrical,Gui_actions}, $L$ has a non-trivial
splitting, in contradiction with property (FA).
This argument extends to the case where $\Gamma$ is a toral relative
hyperbolic group \cite{Groves_limit_hypI}, or where $\Gamma$ has torsion.

\medskip
We end the exposition with the following.
\begin{quest}
\label{que}
Is it true that, at any density, with overwhelming probability all finite index subgroups of a random group satisfy property $\mathrm{(FA)}$?
\end{quest}

If $d>\frac{1}{3}$ then this question has positive answer, since Kazhdan's property (T) implies property (FA) for all finite index subgroups. But already for $d<\frac{1}{5}$, with overwhelming probability a random group does not have property (T) (see \cite[Corollary 7.5]{OWb}). Hence the answer to Question \ref{que} cannot be only based on property (T). If we fix the index of the subgroups considered, Question \ref{que} might have a positive answer justifiable in the spirit of our article. But we expect that the answer to Question \ref{que} in general is much harder.

\bigskip
Our strategy of proof of Theorem \ref{heart} is the following. In the first part we describe a condition which guarantees property $\mathrm{(FA)}$. This part is inspired by an argument of Pride \cite{Pr}, who gives examples of finitely presented groups of cohomological dimension $2$ with property $\mathrm{(FA)}$. (We learned Pride's argument from an article by Delzant and Papasoglu \cite[Theorem 4.1]{DP}.)

More precisely, we prove that certain finite collection of sets of words (the languages of what we call $\frac{1}{3}$--large basic automata) has the property that if we have at least one relator from each of those sets in the presentation of a group, then this group satisfies $\mathrm{(FA)}$. If we compute densities of those sets, they turn out to converge to $1$, if the number of generators in the presentation converges to $\infty$. Hence this argument suffices to prove Corollary \ref{heart_weak}, which says that at each density $d$, if the number of generators is sufficiently large w.r.t. $d$, then a random group with overwhelming probability satisfies $\mathrm{(FA)}$.

In the second part we show that, to some extent, random groups with small number of generators have finite index subgroups which are quotients of random groups with large number of generators. This part of the proof is similar to the argument that random groups in the Gromov density model are quotients of the groups in the triangular model (see \cite[Section I.3.g]{O}). We may then use the fact that property $\mathrm{(FA)}$ is inherited by quotients and by supergroups of finite index. This proves that with overwhelming probability random groups at any density and with any number of generators satisfy $\mathrm{(FA)}$ (Theorem \ref{heart}).

\medskip
If we require in Definition \ref{model} that the relators are cyclically reduced, Theorem \ref{heart} is still valid, although the proof requires small changes. We decided to work mainly in the model in which we allow cyclically non-reduced words, since the proof in this setting is slightly simpler and easier to follow. However, we provide also the proof for the other model.

\medskip
The article is organized as follows. In Section \ref{property} we prove Proposition \ref{many generators}, which provides sufficient conditions for property $\mathrm{(FA)}$ and yields Corollary \ref{heart_weak}, which is a special case of Theorem \ref{heart}. In Section \ref{section breeding} we use Proposition \ref{many generators} to prove Theorem \ref{heart} in full generality. In Section \ref{cyclically reduced} we give a proof of Theorem \ref{heart} in the model allowing only cyclically reduced relators.

\medskip
The third author would like to thank Jacek \'Swi\k{a}tkowski for the introduction into the subject, and the people at the Institut de Math\'ematiques de Toulouse, where this work was carried out, for great atmosphere and hospitality.

\section{Random groups with large number of generators}
\label{property}
In this section we find conditions which guarantee property $\mathrm{(FA)}$ (Proposition \ref{many generators}).

\medskip
We use the following language. An \emph{alphabet} $S$ is a finite set. Let $S^{-1}$ denote the set of the formal inverses to the elements in $S$. Abbreviate $S^\pm=S\cup S^{-1}$. Elements of $S^\pm$ are called \emph{letters}. A \emph{word} over the alphabet $S$ is a sequence of letters.

We fix, for the entire section, an alphabet $S$ and we denote $n=|S|$. Below we define a restricted version of a classical notion of an automaton whose set of states is $\{\emptyset\}\cup S^\pm$.

\begin{defin}
\label{automaton}
A \emph{basic automaton} (shortly a \emph{b-automaton}) \emph{over an alphabet $S$ with transition data $\{\sigma_s\}$} is a pair $(S,\{\sigma_s\})$, where $\{\sigma_s\}_{s\in \{\emptyset\}\cup S^\pm}$ is a family of subsets of $S^\pm$.

The \emph{language} of a b-automaton with transition data $\{\sigma_s\}$ is the set of all (nonempty) words over $S$ beginning with a letter in $\sigma_\emptyset$ and such that for any two consecutive letters $ss'$ we have that $s'\in \sigma_s$.

We say that a b-automaton is \emph{$\lambda$--large}, for some $\lambda\in (0,1)$, if $\sigma_\emptyset\neq \emptyset$ and for each $s\in S^\pm$ we have $|\sigma_s|\geq \lambda 2n$.
\end{defin}


\begin{rem}
\label{count1}
\item[(i)] There are exactly $2^{2n(2n+1)}$ b-automata (over the fixed alphabet $S$ with $n=|S|$).
\item[(ii)] If a b-automaton is \emph{$\lambda$--large}, then its language contains at least $\lceil \lambda 2n \rceil^{L-1}$ words of length $L$ and at least $(\lceil \lambda 2n \rceil-1)^{L-1}$ reduced words of length $L$.
\end{rem}

These estimates are useful in view of the following discussion.

\begin{defin}
\label{growth}
Let $I\subset \mathbf{N}_+$ and let $\mathcal{L}$ be a set of reduced words over an alphabet $S$, containing for all but finitely many $L\in I$ at least $ck^{L}$ words of length $L$, where $c>0, k>1$. Then we say that the \emph{$I$--growth rate of $\mathcal{L}$ is at least $k$}.
\end{defin}

Note that the value $k$ is related to the classical notion of \emph{density} $d_{\mathcal{L}}$ of the set $\mathcal{L}$ by the relation $k=(2n-1)^{d_{\mathcal{L}}}$.
A well known fact in random groups asserts that a random set of relators at density $d$ intersects a fixed set of words of density greater than $1-d$. In the language of the growth rate this fact amounts to the following.

\begin{lemma}[{\cite[Section 9.A]{Gro93}}]
\label{density}
Let $\mathcal{L}$ be a set of reduced words over the alphabet $S$, of the $I$--growth rate at least $k>(2n-1)^{1-d}$. Then with $I$--overwhelming probability, a random set of relators at density $d$ intersects $\mathcal{L}$.
\end{lemma}

We obtain the following corollary. Note that for fixed $d$ and $\lambda$, its hypothesis is satisfied for sufficiently large $n$.

\begin{cor}
\label{big density}
If $\lceil\lambda2n\rceil-1>0$ and $(2n-1)^d\geq\frac{2}{\lambda}$, then with overwhelming probability a random set of relators at density $d$ intersects the languages of all $\lambda$--large b-automata over the alphabet $S$.
\end{cor}

\proof
By Remark \ref{count1}(ii) the $\mathbf{N}_+$--growth rate of the set of reduced words in the language $\mathcal{L}$ of a $\lambda$--large b-automaton is at least $k=\lceil \lambda 2n\rceil -1$. Since
$$2(\lceil \lambda 2n\rceil-1)=(\lceil \lambda 2n\rceil-2)+\lceil \lambda 2n\rceil\geq \lceil \lambda 2n\rceil> \lambda(2n-1),$$ we get that $k> \frac{\lambda}{2}(2n-1)$ which is by hypothesis at least $(2n-1)^{1-d}$. Hence by Lemma \ref{density} a random set of relators at density $d$ intersects $\mathcal{L}$ with overwhelming probability. By Remark \ref{count1}(i), the number of b-automata over the fixed alphabet $S$ depends only on $n$ (and not on $L$), and we get the same conclusion for all languages simultaneously.
\qed

\medskip
Now we present the main result of this section.

\begin{prop}
\label{many generators}
Let $G$ be a group with presentation $\langle S| R\rangle$ such that $R$ intersects the languages of all $\frac{1}{3}$--large b-automata over the alphabet $S$. Then $G$ satisfies $\mathrm{(FA)}$.
\end{prop}

\begin{rem}
\label{constant} The fact that the value $\frac{1}{3}$ of the largeness constant in the hypothesis of Proposition \ref{many generators} is greater than $\frac{1}{4}$ will be surprisingly crucial in the proof of Theorem \ref{heart}. In Section \ref{cyclically reduced}, we will need a modified version of Proposition \ref{many generators}, where the parameter $\frac{1}{3}$ gets closer to $\frac{1}{4}$.
\end{rem}

Before we give the proof of Proposition \ref{many generators}, let us deduce the following consequence, which is a weak version of Theorem \ref{heart}.

\begin{cor}
\label{heart_weak}
Let $0<d<1$ and let $n$ satisfy $(2n-1)^d\geq 6$. Then with overwhelming probability a random group with $n$ generators at density $d$ satisfies property $\mathrm{(FA)}$.
\end{cor}
\proof
Let $\lambda=\frac{1}{3}$. Since $n\geq 2$, we have $\lceil \lambda 2n \rceil -1\geq\lceil \frac{4}{3} \rceil -1>0$. Moreover, $(2n-1)^d\geq 6=\frac{2}{\lambda}$.
By Corollary \ref{big density}, with overwhelming probability $R$ intersects the languages of all $\frac{1}{3}$--large b-automata over the alphabet $S$. Hence by Proposition \ref{many generators} we have that w.o.p. $G$ satisfies $\mathrm{(FA)}$.
\qed

\medskip
The proof of Proposition \ref{many generators} relies on the following lemmas.

\medskip
Note that it is well known that random groups have trivial abelianization, hence they do not admit an epimorphism onto $\mathbf{Z}$. However, it is convenient for us to include the proof of the following version of this assertion.

\begin{lemma}
\label{epimorphism lemma}
If $G=\langle S|R \rangle$ admits an epimorphism onto $\mathbf{Z}$, then there is a $\frac{1}{2}$--large b-automaton over $S$ with language disjoint with $R$.
\end{lemma}

Let us adopt the convention that if
$s\in S^\pm$ (resp. if $w$ is a word over the alphabet $S$), then by $\overline{s}$ (resp. $\overline{w}$) we denote the corresponding element in the group $G$.

\proof
If there is an epimorphism $\psi \colon G \rightarrow \mathbf{Z}$, we consider the following sets. Let
\begin{align*}
S_+&=\{s\in S^\pm \ \mathrm{such} \ \mathrm{that}\  \psi(\overline{s})>0\}, \\
S_0&=\{s\in S^\pm \ \mathrm{such} \ \mathrm{that}\  \psi(\overline{s})=0\}.
\end{align*}
Note that $|S_+\cup S_0|\geq n$ and $S_+$ is nonempty. Consider the b-automaton $\mathds{A}$ over the alphabet $S$ with transition data  $\sigma_\emptyset=S_+$ and $\sigma_s=S_+\cup S_0$, for $s\neq \emptyset$. Let $w=s_{i_1}\ldots s_{i_L}$ be a word in the language of $\mathds{A}$. Then $\psi(\overline{w})=\psi(\overline{s}_{i_1})+\ldots +\psi(\overline{s}_{i_L})>0$. In particular $\overline{w}\neq 0$, hence $w\notin R$. On the other hand, the b-automaton $\mathds{A}$ is $\frac{1}{2}$--large, as required.
\qed

\medskip
Before stating the next lemmas, we need the following discussion of free products with amalgamation.
If $G$ splits as $A\ast_C B$, then we say that an element $g\in G$ is written in a \emph{reduced form} $g=a_1b_1a_2b_2\ldots a_kb_k$ w.r.t. this splitting, if $a_i\in A, b_i\in B$, and none of the terms $a_i,b_i$ belong to $C$, with the exceptions that $a_1,b_k$ are allowed to be trivial in $G$, and that if $g\in C$, we allow $k=1,a_1=g,b_1=e$. The \emph{length} of $g\in G$ w.r.t. the splitting $A\ast_C B$ is the number of the terms $a_i,b_i$ appearing in the reduced form of $g$. The length is well defined, i.e. it does not depend on the reduced form we choose (see \cite[Section I.1.2]{S}). In particular, if $g$ has a reduced form with at least $2$ terms, then it is a nontrivial element of $G$.

If $G=\langle S|R \rangle$ splits as $A\ast_C B$, we denote by $\mathcal{A}, \mathcal{B}, \mathcal{C}, \mathcal{D}$ the sets of letters $s\in S$ whose corresponding $\overline{s}$ lie, respectively, in $A\setminus C,\  B\setminus C,\  C,$ and outside $A\cup B$. (In particular, for $s\in \mathcal{D}$ we have that the length of $\overline{s}$ is at least $2$.) We denote by $\alpha$, $\beta$, $\gamma$, $\delta$ the cardinalities of these sets. We abbreviate $\mathcal{A}^\pm=\mathcal{A}\cup\mathcal{A}^{-1}\subset S^\pm$ and similarly for $\mathcal{B}, \mathcal{C}, \mathcal{D}$.

\begin{lemma}
\label{reduction lemma}
If $G=\langle S|R\rangle$ splits as $A\ast_C B$, then there are $\frac{\alpha+\gamma}{n}$--large and $\frac{\beta+\gamma}{n}$--large b-automata with languages disjoint with $R$.
\end{lemma}

\proof
For each $s'\in S$ we consider the $\frac{\alpha+\gamma}{n}$--large b-automaton $\mathds{A}^{s'}$ over the alphabet $S$ with transition data $\sigma_\emptyset=\{s'\}$ and $\sigma_s=\mathcal{A}^\pm\cup\mathcal{C}^\pm$ for $s\neq \emptyset$. We claim that at least for one $s'\in S$, the language of $\mathds{A}^{s'}$ is disjoint with $R$. Otherwise, for every $s'\in S$ there is a relator $r_{s'}\in R$ contained in the language of $\mathds{A}^{s'}$. Since $\overline{r}_{s'}=1$, we obtain $\overline{s'}\in A$. If this holds for every $s'\in S$, we obtain $G\subset A$, contradiction.
The second construction is analogous.
\qed

\begin{lemma}
\label{main lemma}
Assume that $G=\langle S|R\rangle$ splits as $A\ast_C B$ and the splitting is chosen so  that the sum of the lengths of all generators $\overline{s}\in G$, for $s\in S$, w.r.t this splitting is minimal. Then there is a $\frac{1}{2n}\min \{\delta+\beta, \delta +\alpha\}$--large b-automaton with language disjoint with $R$.
\end{lemma}

We illustrate the idea of the proof by means of the following example. Assume that for all $s\in S$ we have that $\overline{s}\in \mathcal{D}$ and that the first term of the reduced form of $\overline{s}$ lies in $A\setminus C$, and its last term lies lies in $B\setminus C$. (One can check that in this case the minimality hypothesis is satisfied.) Consider the b-automaton with all $\sigma_s$ equal $S\subset S^\pm$. This b-automaton is $\frac{1}{2}$--large (its language consists of all "positive" words). Any word $w=s_{i_1}\ldots s_{i_L}$ in the language of this b-automaton has the following property. If we concatenate reduced forms of all $\overline{s}_{i_l}$, we obtain $\overline{w}$ in a reduced form (there are no cancellations). Thus $\overline{w}$ has large length and cannot be trivial. Hence $w\notin R$. Thus we have constructed a $\frac{1}{2}$--large b-automaton, whose language does not intersect $R$.

\medskip
Before we give the proof of Lemma \ref{main lemma}, we need the following reformulation of the minimality assumption.

\begin{sublemma}
\label{minimality}
Assume that $G=\langle S|R\rangle$ splits as $A\ast_C B$ and the splitting is chosen so that the sum of the lengths of all generators $\overline{s}\in G$, for $s\in S$, w.r.t this splitting is minimal. Then for each $a\in A\setminus C$ we have the following. There are at most $\beta+\delta$ letters $s\in\mathcal{D}^\pm$ with the property that the reduced form of the corresponding $\overline{s}\in G$ begins with a term $a_1\in A\setminus C$, such that $a^{-1}a_1\in C$.
Similarly, for each $b\in B\setminus C$ we have that there are at most $\alpha+\delta$ letters $s\in\mathcal{D}^\pm$ with the property that the reduced form of the corresponding $\overline{s}\in G$ begins with a term $b_1\in B\setminus C$, such that $b^{-1}b_1\in C$.
\end{sublemma}

Note that the first term of the reduced form (which in the notation of the definition of the reduced form is $a_1$, or $b_1$ if $a_1=e$) is determined uniquely modulo multiplying by an element from $C$ on the right (see \cite[Section I.1.2]{S}). We will write shortly \emph{modulo $C$} instead of "modulo multiplying by an element from $C$ on the right".

\proof
We prove the first assertion (the proof of the second one is analogous). Denote by $\mathcal{F}=\mathcal{F}(a)$ the set of all letters in $\mathcal{D}^\pm$, whose reduced forms begin with the term $a$ modulo $C$. Let $\phi=|\mathcal{F}|$. Informally, $\phi$ is the number of "generators' extremities" whose reduced form cancels with $a^{-1}$ or $a$ depending on whether the "extremity" is the "beginning" or the "ending" of the generator.

Let us compute, how do lengths of generators change under conjugating by $a$ (this is equivalent to computing the lengths of the generators with respect to the splitting obtained by conjugating the splitting $A\ast_C B$ by $a$).
For $s\in \mathcal{A}\cup \mathcal{C}$, the length of $a^{-1}\overline{s}a$ equals 1 which is the length of $\overline{s}$. For $s\in \mathcal{B}$, the length of $a^{-1}\overline{s}a$ equals 3, hence increases by $2$ in comparison with the length of $\overline{s}$, which is $1$. For $s\in \mathcal{D}$ we study separately both "extremities" of $\overline{s}$, which means that we study the first terms of reduced forms of $\overline{s}$ for all letters $s\in \mathcal{D}^\pm$. Exactly $\phi$ of these first terms equal $a$ modulo $C$. The other ones are either in $A\setminus (C\cup aC)$ or in $B\setminus C$, and we can denote their numbers, respectively, by $\phi'$ and $2\delta-\phi-\phi'$. This means that conjugating by $a$ increases the sum of the lengths of the generators in $\mathcal{D}$ by $-\phi+(2\delta-\phi-\phi')$.

To summarize, conjugating the splitting $A\ast_C B$ by $a$ gives us a new splitting in which the sum of lengths of generators increases by at most $2\beta-\phi+(2\delta-\phi-\phi')$. By minimality hypothesis on the splitting $A\ast_C B$, we get that this number is non-negative. Since $\phi'\geq 0$, this gives $\phi\leq \beta +\delta$, as required.
\qed

\medskip\par\noindent\textbf{Proof of Lemma \ref{main lemma}.}\ignorespaces
\ We define the following b-automaton $\mathds{A}$ over the alphabet $S$.

Let $\sigma_\emptyset=\mathcal{A}^\pm\cup\mathcal{B}^\pm\cup\mathcal{D}^\pm$.

For $s\in \mathcal{A}^\pm$ let $\sigma_s$ be the union of $\mathcal{B}^\pm$ and those letters $s'\in\mathcal{D}^\pm$ for which the reduced form of $\overline{s'}$ does not begin with $\overline{s}^{-1}$ modulo $C$. Similarly, for $s\in \mathcal{B}^\pm$ let $\sigma_s$ be the union of $\mathcal{A}^\pm$ and those letters $s'\in\mathcal{D}^\pm$ for which the reduced form of $\overline{s'}$ does not begin with $\overline{s}^{-1}$ modulo $C$.

Now suppose that $s\in \mathcal{D}^\pm$, and that $\overline{s}$ ends, in the reduced form, with a term $a_k\in A\setminus C$. Then let $\sigma_s$ be the union of $\mathcal{B}^\pm$ and the set of letters $s'\in\mathcal{D}^\pm$ for which the reduced form of $\overline{s'}$ does not begin with $a_k^{-1}$ modulo $C$. Analogously, if $s\in \mathcal{D}^\pm$, and $\overline{s}$ ends, in the reduced form, with a term $b_k\in B\setminus C$, then let $\sigma_s$ be the union of $\mathcal{A}^\pm$ and the set of letters $s'\in\mathcal{D}^\pm$ for which the reduced form of $\overline{s'}$ does not begin with $b_k^{-1}$ modulo $C$.

Finally, for $s\in \mathcal{C}^\pm$ let $\sigma_s=S^\pm$.

\medskip\par\noindent\textbf{Step 1.}\ignorespaces
\
\emph{The b-automaton $\mathds{A}$ is $\frac{1}{2n}\min \{\delta+\beta, \delta +\alpha\}$--large.}

\medskip
If  $s\in \mathcal{A}^\pm$, then, by Sublemma \ref{minimality}, we have in $\mathcal{D}^\pm$ at least $2\delta - (\beta+\delta)$ elements of $\sigma_s$. Adding elements of $\mathcal{B}^\pm$, we get altogether at least $\delta+\beta$ elements of $\sigma_s$.
Analogously, if $s\in \mathcal{B}^\pm$, there are at least $\delta+\alpha$ elements in $\sigma_s$.
The computation is similar for $s\in\mathcal{D}^\pm$.

Cases where $s\in \mathcal{C}^\pm$ or $s=\emptyset$ are obvious.

\medskip\par\noindent\textbf{Step 2.}\ignorespaces
\ \emph{The language of $\mathds{A}$ is disjoint with $R$.}

\medskip
It is enough to prove that for any word $w$ in the language of $\mathds{A}$, the element $\overline{w}\in G$ is nontrivial. This follows from the following stronger assertion.

\medskip\par\noindent\textbf{Claim.}\ignorespaces
\ \emph{For any word $w$ of length $k$ in the language of $\mathds{A}$, the length of $\overline{w}$ w.r.t. the splitting $A\ast_C B$ is at least $k$ and the following holds. If we denote by $s$ the last letter of $w$, we have that the last term of the reduced form of $\overline{w}$ equals the last term of the reduced form of $\overline{s}$ modulo multiplying by an element from $C$ on the left.}

\medskip
We prove the claim by induction on $k$. For $k=1$ this is obvious. Assume we have already proved the claim for $k=l-1\geq 1$. Now let $w$ be a word of length $l$ in the language of $\mathds{A}$ ending with the pair $s's$. Then the word $w'$ obtained from $w$ by removing $s$ from the end also lies in the language of $\mathds{A}$ and we can apply to it the induction hypothesis. We get that $\overline{w'}$ has length at least $l-1$ and its reduced form ends with the last term, say $b\in B\setminus C$, of the reduced form of $\overline{s'}$.

If $s\in \mathcal{D}^\pm$, then the length of $\overline{s}$ is at least $2$ and by definition of $\sigma_{s'}$ the first term $t$ of the reduced form of $\overline{s}$ does not equal $b^{-1}$ modulo multiplying by an element from $C$ on the right. Hence when we concatenate the reduced forms of $\overline{w'}$ and $\overline{s}$ and, if $t\in B\setminus C$, when we substitute the pair $bt$ with a single term in $B\setminus C$, we obtain $\overline{w}$ in a reduced form of length at least $l$ and whose last term equals the last term of the reduced form of $\overline{s}$, as required.

By definition of $\sigma_{s'}$, the only other possibility for $s$ is that it lies in $\mathcal{A}^\pm$, i.e. that $\overline{s}\in A\setminus C$. Thus adjoining $\overline{s}$ at the end of the reduced form of $\overline{w'}$ gives a reduced form of $\overline{w}$, which is of length at least $l$ and whose last term equals the last (and only) term of the reduced form of $\overline{s}$. This proves the claim for $k=l$, and ends the induction proof.

\medskip
This ends the proof of Lemma \ref{main lemma}.
\qed

\medskip
We now collect all pieces of information.

\medskip\par\noindent\textbf{Proof of Proposition \ref{many generators}.}\ignorespaces
\ Since $G$ is finitely generated, we need to prove that $G$ does not admit an epimorphism onto $\mathbf{Z}$ and does not split as a free product with amalgamation.
By Lemma \ref{epimorphism lemma}, since $\frac{1}{2}> \frac{1}{3}$, we have that $G$ does not admit an epimorphism onto $\mathbf{Z}$. It remains to prove that $G$ does not split as a free product with amalgamation. We prove this by contradiction.

Assume that $G$ splits as $A\ast_C B$ and the splitting is chosen so that the sum of the lengths of all generators $\overline{s}\in G$, for $s\in S$, w.r.t this splitting is minimal.
By Lemma \ref{main lemma} we have that $\frac{1}{2n}\min \{\delta+\beta, \delta +\alpha\}<\frac{1}{3}$.
Assume, w.l.o.g., that $\beta\leq \alpha$. Then $\delta+\beta<\frac{2n}{3}$.

By Lemma \ref{reduction lemma} we have that $\frac{\alpha+\gamma}{n}<\frac{1}{3}$, i.e. that $\alpha+\gamma<\frac{n}{3}$.
Adding up, we obtain $\delta+\beta+\alpha+\gamma<n$, contradiction.
\qed

\section{Increasing the number of generators}
\label{section breeding}
In this section we demonstrate how to pass from a model where random groups have small number of generators to a model with large number of generators, where we can apply Proposition \ref{many generators}.

\medskip Recall that in our random model we denote the set of generators by $S$ with $n=|S|$. The density of the random set of relators is denoted by $d$. For our argument we need to fix some natural number $B$ which we will later require to be sufficiently large so that $\sqrt[B]{12}<(2n-1)^d$ (this estimate will be used only once at the end of the proof).

Let $\widetilde{S}$ denote the set of reduced words of length $B$ over the alphabet $S$.
The involution on $\widetilde{S}$ mapping each word to its inverse does not have fixed points. Thus we can partition $\widetilde{S}$ into $\hat{S}$ and $\hat{S}^{-1}$. We denote  $\hat{S}^\pm=\hat{S}\cup\hat{S}^{-1}$ (instead of $\widetilde{S}$).
Let $\hat{n}$ be the cardinality of $\hat{S}$, which equals $n(2n-1)^{B-1}$.

Recall that $L$ denotes the length of the random relators. Our proof is significantly simpler, if we consider only those $L$ that are divisible by $B$. We will always distinguish this case and we recommend the reader to focus on this case during the first reading of the article. For $0\leq P<B$ let $I_P\subset \mathbf{N}_+$ denote the set of those $L$ that can be written as $L=B\hat{L}+P$.


\begin{defin}
\ Let $r$ be a word of length $L\in I_0$ over the alphabet $S$. Divide the word $r$ into $\hat{L}$ blocks of length $B$. This determines a new word $\hat{r}$ of length $\hat{L}$ over the alphabet $\hat{S}$, which we call the word \emph{associated} to $r$.
\end{defin}

\begin{defin}
Given a set $R$ of relators over $S$ of equal length $L\in I_P$, we define the \emph{associated group $\hat{G}$} in the following way.

If $P=0$, then we consider the set $\hat{R}$ of relators associated to relators in $R$. We define $\hat{G}$ to be the group $\langle \hat{S}| \hat{R}\rangle$.

If $1\leq P<B$, then there is no natural way to associate relators over $\hat{S}$ to relators over $S$. We resolve this in the following way. Suppose that $r_1,r_2\in R$ are two relators of length $L$ over $S$, such that $r_1=q_1v^{-1}$ and $r_2=vq_2$, for some word $v$ over $S$ of length $P$. We then obtain a (possibly non-reduced) word $q_1q_2$ over $S$, of length $2B\hat{L}$, with the property that $\overline{q_1}\overline{q_2}=1$ in $G=\langle S|R\rangle$. To this word we can associate, as before, a relator over $\hat{S}$, of length $2\hat{L}$, which we denote by $\hat{r}(r_1,r_2)$.
We denote by $\hat{R}$ the set of all $\hat{r}(r_1,r_2)$ as above and we define $\hat{G}=\langle \hat{S}| \hat{R}\rangle$.
\end{defin}

\begin{lemma}
\label{inheriting FA}
If $\hat{G}$ satisfies property $\mathrm{(FA)}$, then so does $G=\langle S| R\rangle$.
\end{lemma}

\proof
Indeed, we have a natural epimorphism $\hat{G}\rightarrow H$, where $H$ is the subgroup of $G$ generated by the words of length $B$ over $S$. If $\hat{G}$ satisfies $\mathrm{(FA)}$, then its quotient $H$ also satisfies $\mathrm{(FA)}$. Moreover, since $H\subset G$ is of finite (in fact at most $2n$) index, we have by \cite[Section I.6.3.4]{S} that $G$ also satisfies $\mathrm{(FA)}$.
\qed

\medskip
The idea behind the remaining part of the proof is the following. Assume that $L=B\hat{L}$. Then each relator of length $\hat{L}$ over the alphabet $\hat{S}$ is associated to a relator of length $L$ over the alphabet $S$. Consider the case of a model where we allow non-reduced relators. Then one can check that the density of a set $R$ of relators over $S$ equals the density of the set $\hat{R}$ of associated relators over $\hat{S}$. This means that $\hat{G}$ is a random group at density $d$ in a model, where we allow non-reduced relators, with a large number $\hat{n}=|\hat{S}|$ of generators.
Hence in this context, Corollary \ref{heart_weak} implies, in view of Lemma \ref{inheriting FA}, Theorem \ref{heart}. However, we have decided not to work in this model, since it is less standard. For some reference on it, see \cite{Oll04}.

In our setting, we have to resolve the problem that some reduced words over $\hat{S}$ might be associated to non-reduced words over $S$. The key is the following.


\begin{lemma}
\label{counting}
Let $\mathds{A}$ be a $\lambda$--large b-automaton over $\hat{S}$. Denote by $\mathcal{L}_\mathds{A}$
the set of reduced words over the alphabet $S$, whose associated words lie in the language of $\mathds{A}$. Assume that $\lambda'=\lambda-\frac{1}{2n}>0$. Then the $I_0$--growth rate of $\mathcal{L}_\mathds{A}$ is at least $\sqrt[B]{\lambda'}(2n-1)$.
\end{lemma}

The outline of the proof is the following. We construct a b-automaton $\mathds{A}^\red$ over $\hat{S}$, whose language consists of
elements of the language of $\mathds{A}$, which are associated to reduced words over $S$. In other words, the language of $\mathds{A}^\red$ consists of words associated to elements of $\mathcal{L}_{\mathds{A}}$. Then we estimate from below the growth rate of the language of $\mathds{A}^\red$, hence the growth rate of $\mathcal{L}_{\mathds{A}}$, in terms of $n, B,$ and $\lambda$.

\proof Denote the transition data of $\mathds{A}$ by $\{\sigma_{\hat{s}}\}$. For any $\hat{s}\in \hat{S}^\pm$, let $\rho_{\hat{s}}\subset \hat{S}^\pm$ be the set of
$\hat{s}'$ such that the first letter of $\hat{s}'$ interpreted as a word over the alphabet $S$ is the inverse of the last letter of $\hat{s}$ as a word over $S$. Observe that $|\rho_{\hat{s}}|=\frac{1}{2n}2\hat{n}$.

Let $\mathds{A}^{\red}$ be the b-automaton over the alphabet $\hat{S}$ with transition data $\sigma^{\red}_\emptyset=\sigma_\emptyset$ and $\sigma^{\red}_{\hat{s}}=\sigma_{\hat{s}}\setminus \rho_{\hat{s}}$, for $\hat{s}\in \hat{S}^\pm$. We have that $|\sigma^\red_{\hat{s}}|\geq |\sigma_{\hat{s}}|-|\rho_{\hat{s}}|\geq \lambda 2\hat{n}-\frac{1}{2n}2\hat{n}=\lambda' 2\hat{n}$. Hence $\mathds{A}^\red$ is $\lambda'$--large.
By Remark \ref{count1}(ii), its language contains at least $(\lceil \lambda'2\hat{n} \rceil)^{\hat{L}-1}$ words of length $\hat{L}$.

\medskip
Observe that the language of $\mathds{A}^\red$ has the following two properties. First, it is contained in the language of $\mathds{A}$. Second, for any word $w$ of length $\hat{L}$ in the language of $\mathds{A}^\red$, if we substitute each letter $\hat{s}\in w$ with the corresponding word over the alphabet $S$, we obtain a reduced word of length $B\hat{L}$ over $S$.

This implies that the number of reduced words of length $L=B\hat{L}$ over the alphabet $S$, whose associated words lie in the language of $\mathds{A}$ is bounded from below by
$$\left(\left\lceil \lambda'2\hat{n}\right\rceil\right)^{\hat{L}-1}\geq
\left(\lambda'2\hat{n}\right)^{\frac{L}{B}-1}
\geq c\left(\sqrt[B]{\lambda'2\hat{n}}\right)^L,$$
for some $c>0$. In other words, the $I_0$--growth rate of $\mathcal{L}_{\mathds{A}}$ is at least $k=\sqrt[B]{\lambda'2\hat{n}}$. By definition of $\hat{n}$ we have that
$$k>\sqrt[B]{\lambda'(2n-1)^B}=\sqrt[B]{\lambda'}(2n-1).$$
This ends the proof of Lemma \ref{counting}.
\qed

\medskip
We obtain some corollaries for the case where $P\neq 0$. We recommend the reader only interested in the case of $L$ divisible by $B$
to skip them and proceed directly to the proof of Theorem \ref{heart}.

\medskip
Assume that $1\leq P<B$. We keep the setting from Lemma \ref{counting}. Let $\mathcal{P}^P_{\mathds{A}}$ (the "prefix set") denote the set of reduced words $w$ over the alphabet $S$, whose length $L$ lies in $I_P$, such that the length $B\hat{L}=L-P$ prefix of $w$ lies in $\mathcal{L}_{\mathds{A}}$.

\begin{cor}
\label{prefix}
The $I_P$--growth rate of $\mathcal{P}^P_{\mathds{A}}$ is at least $k>\sqrt[B]{\lambda'}(2n-1)$.
\end{cor}
\proof
This follows from the fact that any word of length $B\hat{L}$ in $\mathcal{L}_{\mathds{A}}$ can be extended to a word of length $L$ in $\mathcal{P}^P_{\mathds{A}}$, and from Lemma \ref{counting}.
\qed

\medskip
Consider the set of reduced words of length $B$ over the alphabet $S$, which begin with $s^{-1}$, for some $s\in S^\pm$. View this set as a subset $\rho^s$ of $\hat{S}^\pm$. For any $\hat{s}\in \hat{S}^\pm$ let $\mathds{A}^{\hat{s},s}$ be the b-automaton over the alphabet $\hat{S}$ with transition data equal to the transition data of $\mathds{A}$ with the exception that we substitute $\sigma_{\emptyset}$ with $\sigma_{\hat{s}}\setminus \rho^s$. This set is nonempty since $|\sigma_{\hat{s}}|\geq \lambda 2\hat{n}$ and $|\rho^s|=\frac{1}{2n}2\hat{n}$. Hence $\mathds{A}^{\hat{s},s}$ is $\lambda$--large.

Let $v$ be a reduced word of length $P$ over the alphabet $S$ ending with a letter $s$.
Let $\mathcal{S}_{\mathds{A}}^{\hat{s},v}$ (the "suffix set") denote the set of words $w$ over the alphabet $S$, whose length $L$ lies in $I_P$, such that the length $P$ prefix of $w$ equals $v$, and the length $B\hat{L}=L-P$ suffix of $w$ lies in $\mathcal{L}_{\mathds{A}^{\hat{s},s}}$.
Since all words in the language of $\mathds{A}^{\hat{s},s}$ start with a letter outside $\rho^s$, we have that all words in $\mathcal{L}_{\mathds{A}^{\hat{s},s}}$ start with a letter different from $s^{-1}$, and consequently all words in $\mathcal{S}_{\mathds{A}}^{\hat{s},v}$ are reduced.

By Lemma \ref{counting} applied to $\mathds{A}^{\hat{s},s}$ we obtain immediately the following.

\begin{cor}
\label{sufffix}
The $I_P$--growth rate of $\mathcal{S}_{\mathds{A}}^{\hat{s},v}$ is at least $k>\sqrt[B]{\lambda'}(2n-1)$.
\end{cor}

We are now ready for the following.

\medskip\par\noindent\textbf{Proof of Theorem \ref{heart}.}\ignorespaces
\
Let $R$ denote a random set of relators over $S$ and $G=\langle S|R\rangle$. We choose $B$ sufficiently large so that $\sqrt[B]{12}<(2n-1)^d$. Let $\hat{G}=\langle \hat{S}|\hat{R}\rangle$ be the associated group.
We want to verify, with overwhelming probability, the hypothesis of Proposition \ref{many generators} for $\hat{G}$, that the set $\hat{R}$ intersects the languages of all $\frac{1}{3}$--large b-automata over the alphabet $\hat{S}$. By Remark \ref{count1}(i) it is enough to prove this for a single $\frac{1}{3}$--large b-automaton $\mathds{A}$. As above, we denote by $\mathcal{L}_\mathds{A}$ the set of reduced words over the alphabet $S$, whose associated words lie in the language of $\mathds{A}$. By Lemma \ref{counting} the $I_0$--growth rate of $\mathcal{L}_\mathds{A}$ is at least $k>\sqrt[B]{\lambda'}(2n-1)$, where (since $n\geq 2$) we have
$$\lambda'=\frac{1}{3}-\frac{1}{2n}\geq \frac{1}{3}-\frac{1}{4}=\frac{1}{12}.$$ (This is the point to which we refer in Remark \ref{constant}.)

By the choice of $B$, the $I_0$--growth rate of $\mathcal{L}_\mathds{A}$ is at least $\frac{2n-1}{\sqrt[B]{12}}>(2n-1)^{1-d}$.
Thus we can apply Lemma \ref{density} and we get that
with $I_0$--overwhelming probability there is a relator $r\in R\cap \mathcal{L}_{\mathds{A}}$, hence there is a relator $\hat{r}\in \hat{R}$ in the language of $\mathds{A}$, as required.

\medskip
We have thus proved that the hypothesis of Proposition \ref{many generators} for the group $\hat{G}$ is satisfied with $I_0$--overwhelming probability. In that case, by Proposition \ref{many generators},  $\hat{G}$ satisfies property $\mathrm{(FA)}$. By Lemma \ref{inheriting FA} this implies that $G$ satisfies $\mathrm{(FA)}$. This ends the proof of Theorem \ref{heart} under the assumption that we consider only $L\in I_0$.

\medskip
We now focus on the remaining case where $L\in I_P$ with $P\neq 0$. Since the number of the sets $\mathcal{S}_{\mathds{A}}^{\hat{s},v}$ (defined before Corollary \ref{sufffix}) is finite and independent of $L$, we get by Corollaries \ref{prefix} and \ref{sufffix}, by the choice of $B$, and by Lemma \ref{density}, that with $I_P$--overwhelming probability a random set of relators $R$ contains an element in $\mathcal{P}^P_{\mathds{A}}$ and elements in $\mathcal{S}_{\mathds{A}}^{\hat{s},v}$, for all $\hat{s},v$.

In that case let $r_1\in R\cap\mathcal{P}^P_{\mathds{A}}$. Denote by $v^{-1}$ the word consisting of last $P$ letters of $r_1$ and by $\hat{s}$ the letter in $\hat{S}^\pm$ associated to the length $B$ block appearing before $v^{-1}$ in $r_1$. Let $r_2\in R\cap \mathcal{S}_{\mathds{A}}^{\hat{s},v}$. Then the relator $\hat{r}(r_1,r_2)$ belongs to both $\hat{R}$ and the language of the b-automaton $\mathds{A}$.

\medskip
Hence, with $I_P$--overwhelming probability, the hypothesis of Proposition \ref{many generators} is satisfied and we can conclude that $\hat{G}$ satisfies property $\mathrm{(FA)}$. By Lemma \ref{inheriting FA} this implies property $\mathrm{(FA)}$ for $G$ and ends the proof of Theorem \ref{heart} in the case $L\in I_P$ for $P\neq 0$.
\qed

\section{Cyclically reduced relators model}
\label{cyclically reduced}
In this section we explain what changes need to be introduced in the proof of Theorem \ref{heart} in the case where we require random relators to be cyclically reduced.

\begin{theorem}
\label{heart_reduced}
Let $0<d<1$. Then in the model in which we allow only cyclically reduced relators, with overwhelming probability a random group at density $d$ satisfies property $\mathrm{(FA)}$.
\end{theorem}

The problem is that the words in $\mathcal{L}_{\mathds{A}}$ (see the proof of Theorem \ref{heart}) might be not cyclically reduced. Moreover, we do not have a guarantee that for a given word of length $B(\hat{L}-1)$ in $\mathcal{L}_{\mathds{A}}$, we can extend this word to any word of length $B\hat{L}$ in $\mathcal{L}_{\mathds{A}}$, which is cyclically reduced. This spoils the counting in Lemma \ref{counting}.

To overcome this, we need to consider slightly wider class of automata than we have used so far,
with richer languages.
Shortly, we allow a different transition rule for the last letter, a rule that allows almost half of the letters to be put on the end.

\begin{defin}
An \emph{enhanced basic automaton} (shortly an \emph{e-automaton}) over an alphabet $S$ is a b-automaton over $S$ together with a \emph{final transition data} $\{\tau_s\}$, which is a family of sets $\tau_s\subset S^\pm$ for all $s\in S^\pm$.

The \emph{language} of an e-automaton is the set of all (nonempty) words in $S$ beginning with a letter in $\sigma_\emptyset$ and such that for any two consecutive letters $ss'$ we have that $s'\in \sigma_s$, if $s'$ is not the last letter, and $s'\in \tau_s$, if $s'$ is the last letter.

We say that an e-automaton is \emph{$(\lambda,\e)$--large}, for some $\lambda\in (0,1), \e \in (0,\frac{1}{2})$, if its underlying b-automaton is \emph{$\lambda$--large}, and for each $s\in S^\pm$ we have $|\tau_s|> (\frac{1}{2}-\e)2n$
(the reason the latter condition is expressed in this way will become clear in Sublemma \ref{epsilon}).
\end{defin}

\begin{rem}
\label{promoting}
If $\lambda>\frac{1}{2}-\e$, then a $\lambda$--large b-automaton can be promoted to a $(\lambda,\e)$--large e-automaton with the same language by just putting $\tau_s=\sigma_s$, for all $s\in S^\pm$.
\end{rem}

First we show that using this notion we can save the counting argument from Lemma \ref{counting}. For an e-automaton $\mathds{A}^e$ over the alphabet $\hat{S}$, we denote by $\mathcal{L}^{\cyc}_{\mathds{A}^e}$ the set of cyclically reduced words over $S$, whose associated words lie in the language of $\mathds{A}^e$. We have the following analogue (and consequence) of Lemma \ref{counting}.

\begin{lemma}
\label{last lemma}
There is a constant $\e_0\in (0,\frac{1}{2})$ such that for any $\e\in (0,\e_0]$ we have the following. Let $B\geq 3$. Let $\mathds{A}^e$ be a $(\lambda,\e)$--large e-automaton with $\lambda'=\lambda -\frac{1}{2n}>0$. Then the $I_0$--growth rate of $\mathcal{L}^{\cyc}_{\mathds{A}^e}$ is at least $\sqrt[B]{\lambda'}(2n-1)$.
\end{lemma}
\proof
Let $\mathds{A}$ be the underlying automaton of $\mathds{A}^e$. By Lemma \ref{counting}, the $I_0$--growth rate of $\mathcal{L}_{\mathds{A}}$ is at least $\sqrt[B]{\lambda'}(2n-1)$.
Hence, to prove Lemma \ref{last lemma} it suffices to show that any length $B(\hat{L}-1)$ word in $\mathcal{L}_{\mathds{A}}$ can be extended to a length $B\hat{L}$ word in $\mathcal{L}^{\cyc}_{\mathds{A}^e}$. We need to compute the following.
\begin{sublemma}
\label{epsilon}
For any alphabet $S$ with $n=|S|\geq 2$, there is a constant $\e_0\in (0,\frac{1}{2})$ such that for any $B\geq 3$ and any letters $s,s'\in S^\pm$, we have the following. Let $2\hat{n}$ be the number of all reduced words of length $B$ in $S$ and let $\hat{n}'$ be the number of reduced words of length $B$ over $S$ which begin with $s$ or end with $s'$. Then $\hat{n}'\leq (\frac{1}{2}-\e_0)2\hat{n}$.
\end{sublemma}

Actually, we can take the same $\e_0=\frac{1}{18}$ for any $n$.

\proof
The number of reduced words which begin with $s$ and the number of reduced words which end with $s'$ equal $\frac{1}{2n}2\hat{n}$, so $\hat{n}'\leq \frac{1}{n}2\hat{n}$. Hence if $n\geq 3$, one can take $\e_0=\frac{1}{2}-\frac{1}{n}$. If $n=2$, we need to estimate the number of words which simultaneously begin with $s$ and end with $s'$. We have at least $(2n-1)^{B-3}(2n-2)=\frac{\hat{n}(2n-2)}{n(2n-1)^2}$ such words. Hence
$$\frac{\hat{n}'}{2\hat{n}} \leq \frac{1}{n}-\frac{(2n-2)}{2n(2n-1)^2}=\frac{1}{2}-\frac{1}{18},$$
and we can take $\e_0=\frac{1}{18}$.
\qed

\medskip
We return to the proof of Lemma \ref{last lemma}. Let $w$ be a length $B(\hat{L}-1)$ word in $\mathcal{L}_{\mathds{A}}$, and let $\hat{w}$ be the word associated to it. Denote the last letter of $\hat{w}$ by $\hat{s}$. Let $s,s'$ denote the first and the last letter of $w$. The number of reduced words of length $B$ over the alphabet $S$ which start with $(s')^{-1}$ or end with $s^{-1}$ is at most $(\frac{1}{2}-\e)2\hat{n}$ by Sublemma \ref{epsilon}. Hence, by definition of largeness, there is a letter in $\tau_{\hat{s}}$ (associated to some word $v$ over $S$) without this property. Thus the word $wv$ lies in
$\mathcal{L}^{\cyc}_{\mathds{A}^e}$. This ends the proof of Lemma \ref{last lemma}.
\qed

\medskip
Now for $1\leq P<B$ let $\mathcal{P}^P_{\mathds{A}}$ be defined as in Section \ref{section breeding} and let $(\mathcal{P}^P)^{\cyc}_{\mathds{A}}\subset \mathcal{P}^P_{\mathds{A}}$ be the subset of cyclically reduced words.
By Lemma \ref{counting} we obviously have the following.

\begin{cor}
\label{prefix2}
The $I_P$--growth rate of $(\mathcal{P}^P)^{\cyc}_{\mathds{A}}$ is at least $k>\sqrt[B]{\lambda'}(2n-1)$.
\end{cor}

Let $\mathcal{S}^{\hat{s},v}_{\mathds{A}^e}$ be defined from $\mathds{A}^e$ in the same way we defined $\mathcal{S}^{\hat{s},v}_{\mathds{A}}$ from $\mathds{A}$ in Section \ref{section breeding}. Let $(\mathcal{S}^{\hat{s},v})^{\cyc}_{\mathds{A}^e}\subset \mathcal{S}^{\hat{s},v}_{\mathds{A}^e}$ be the subset of cyclically reduced words. The following analogue of Corollary \ref{sufffix} can be obtained from Sublemma \ref{epsilon} in the same way as Lemma \ref{last lemma}.

\begin{cor}
\label{sufffix2}
The $I_P$--growth rate of $(\mathcal{S}^{\hat{s},v})^{\cyc}_{\mathds{A}^e}$ is at least $k>\sqrt[B]{\lambda'}(2n-1)$.
\end{cor}

Now we prove the following version of Proposition \ref{many generators}. Let $\e=\min \{ \e_0, \frac{1}{6}\}$, where $\e_0\in (0,\frac{1}{2})$ is the constant from Lemma \ref{last lemma} and Sublemma \ref{epsilon}.

\begin{prop}
\label{many generators_reduced}
Let $G$ be a group with presentation $\langle S| R\rangle$ such that the set of relators in $R$ of length at least $3$ intersects the languages of all $(\frac{1}{4}+\frac{\e}{2},\e)$--large e-automata over the alphabet $S$. Then $G$ satisfies $\mathrm{(FA)}$.
\end{prop}

Before we give the proof, we need the following analogue of Lemma \ref{main lemma}.

\begin{lemma}
\label{main_lemma_reduced}
Assume that $G$ splits as $A\ast_C B$ and the splitting is chosen so that the sum of the lengths of all generators $\overline{s}\in G$, for $s\in S$, w.r.t this splitting is minimal.
Then there is a $(\frac{1}{2n}\min \{\delta+\beta, \delta +\alpha\},\e)$--large e-automaton with language disjoint with the set of relators in $R$ of length at least $3$.
\end{lemma}

\proof
We take the b-automaton $\mathds{A}$ described in the proof of Lemma \ref{main lemma}, and promote it to an e-automaton $\mathds{A}^e$ by putting $\tau_s=\sigma_s\cup \mathcal{A}^\pm\cup\mathcal{B}^\pm\cup \mathcal{C}^\pm$.

\medskip\par\noindent\textbf{Step 1.}\ignorespaces
\ \emph{$\mathds{A}^e$ is $(\frac{1}{2n}\min \{\delta+\beta, \delta +\alpha\},\e)$--large.}

\medskip
By step 1 in the proof of Lemma \ref{main lemma}, we just need to estimate $|\tau_s|$. By Sublemma \ref{minimality} we have that $|\sigma_s\cap\mathcal{D}^\pm|\geq \delta-\beta$ or $|\sigma_s\cap\mathcal{D}^\pm|\geq \delta-\alpha$. Since $$|\tau_s|=|\sigma_s\cap\mathcal{D}^\pm|+|\mathcal{A}^\pm\cup\mathcal{B}^\pm\cup \mathcal{C}^\pm|,$$ we have that $$|\tau_s|\geq \delta-\max \{\alpha, \beta \} +2(\alpha+\beta+\gamma)\geq \alpha+\beta+\gamma +\delta=n>\left(\frac{1}{2}-\e\right)2n,$$ as required.

\medskip\par\noindent\textbf{Step 2.}\ignorespaces
\ \emph{The language of $\mathds{A}^e$ is disjoint with the set of relators in $R$ of length at least $3$.}

\medskip
Otherwise, a word $r\in R$ which is in the language of $\mathds{A}^e$ either lies in the language of the underlying b-automaton $\mathds{A}$ (which is not possible by step 2 in the proof of Lemma \ref{main lemma}) or is a concatenation of a word in the language of $\mathds{A}$ and a letter in $\mathcal{A}^\pm\cup\mathcal{B}^\pm\cup \mathcal{C}^\pm$.
In the latter case, $r$ is of the form $ws$, where $\overline{w}$ has length at least $2$ (by the claim in the proof of Lemma \ref{main lemma}) and $(\overline{s})^{-1}$ has length $1$. Contradiction.
\qed

\medskip
Now we are ready for the following.

\medskip\par\noindent\textbf{Proof of Proposition \ref{many generators_reduced}.}\ignorespaces
\
We argue as in the proof of Proposition \ref{many generators}.
Again we need to prove that $G$ does not admit an epimorphism onto $\mathbf{Z}$ and does not split as a free product with amalgamation.

\medskip
By Lemma \ref{epimorphism lemma}, in view of Remark \ref{promoting} (note that $\frac{1}{2}>\frac{1}{2}-\e$), and since  $\frac{1}{2}>\frac{1}{4}+\frac{\e}{2}$ we have that $G$ does not admit an epimorphism onto $\mathbf{Z}$. It remains to prove that $G$ does not split as a free product with amalgamation. We prove this by contradiction.

\medskip
Assume that $G$ splits as $A\ast_C B$ and the splitting is chosen so that the sum of the lengths of all generators $\overline{s}\in G$, for $s\in S$, w.r.t this splitting is minimal.

By Lemma \ref{main_lemma_reduced} we obtain that $\frac{1}{2n}\min \{\delta+\beta, \delta +\alpha\}<\frac{1}{4}+\frac{\e}{2}$. Assume, w.l.o.g., that $\beta\leq \alpha$. Then $\frac{\delta+\beta}{n}<\frac{1}{2}+\e$, hence $$\frac{\alpha+\gamma}{n}=1-\frac{\delta+\beta}{n}>\frac{1}{2}-\e.$$

By Lemma \ref{reduction lemma} and Remark \ref{promoting}, there is a $(\frac{1}{2}-\e,\e)$--large e-automaton whose language does not intersect $R$. Hence $\frac{1}{2}-\e<\frac{1}{4}+\frac{\e}{2}$, which is equivalent to $\e>\frac{1}{6}$. This contradicts the choice of $\e$.
\qed

\medskip
We conclude with the following.

\medskip\par\noindent\textbf{Proof of Theorem \ref{heart_reduced}.}\ignorespaces
\ Choose $B$ such that $\sqrt[B]{\frac{2}{\e}}<(2n-1)^d$, where $\e$ is the constant from Proposition \ref{many generators_reduced}.
By Lemma \ref{last lemma} applied to $\lambda=\frac{1}{4}+\frac{\e}{2}$,
for any $(\frac{1}{4}+\frac{\e}{2},\e)$--large e-automaton $\mathds{A}$,
the $I_0$--growth rate of $\mathcal{L}^{\cyc}_{\mathds{A}^e}$ is at least $k=\sqrt[B]{\lambda'}(2n-1)$, where $\lambda'=\lambda -\frac{1}{2n}\geq \frac{\e}{2}$.

Since we chose $B$ such that $\frac{1}{\sqrt[B]{\lambda'}}<(2n-1)^d$, we have that this $I_0$--growth rate satisfies $k>(2n-1)^{1-d}$.
By Lemma \ref{density} (which is valid also in this model), w.$I_0$--o.p. the set $\hat{R}$, of relators associated to a random set $R$ of cyclically reduced relators, intersects the languages of all $(\frac{1}{4}+\frac{\e}{2})$--large e-automata over the alphabet $\hat{S}$. Then, by Proposition \ref{many generators_reduced},
the associated group $\hat{G}=\langle \hat{S}| \hat{R}\rangle$ satisfies $\mathrm{(FA)}$ and thus, by Lemma \ref{inheriting FA}, we have that $G=\langle S|R\rangle$ satisfies $\mathrm{(FA)}$.

\medskip
It remains to consider the case of $L\in I_P$, for $P\neq 0$. This case follows from Corollaries \ref{prefix2} and \ref{sufffix2}. The argument is analogous to the one in the proof of Theorem \ref{heart}, where we use Corollaries \ref{prefix} and \ref{sufffix}, and we omit it.
This ends the proof of Theorem \ref{heart_reduced}.
\qed

\end{document}